\documentclass[11pt]{article}
\usepackage{amssymb}
\usepackage{amsmath}
\usepackage{amsthm}
\usepackage{hyperref}

\newtheorem{Theorem}{Theorem}
\newtheorem{Proposition}{Proposition}

\newcommand{\pr}{\mbox{pr}_1\,}

\newcommand{\op}[1]{\mathop{\oplus}\limits_{\phantom{.}#1}}
\newcommand{\opp}[2]{\mathop{\oplus}\limits_{\phantom{.}#1}^{\phantom{.}#2}}

\title{Projection representable relations\\ on Menger $(2,n)$-semigroups
\thanks{\hspace*{5mm}Keywords: $n$-place function, algebra
of functions, Menger algebra, $(2,n)$-semigroup.}}

\author{\sc Wies{\l}aw A. Dudek
and Valentin S. Trokhimenko}

\begin{document}
\sloppy \maketitle

\begin{abstract}\noindent In this
paper relations of non-empty intersection, inclusion end equality
of domains of functions for $(2,n)$-semigroups of partial
$n$-place functions are investigated.
\end{abstract}

\section{Introduction}
Investigation of partial multiplace functions by algebraic methods
plays an important role in modern mathematics where we consider
various operations on sets of functions which are naturally
defined. The basic operation for $n$-place functions is a
superposition (composition) $O$ of $n+1$ such functions, but there
are some other naturally defined operations, which are also worth
considering. In this paper we consider binary Mann's compositions
$\op{1},\ldots,\op{n}$ for partial $n$-place functions introduced
in \cite{Man}, which have many important applications for the
studies of binary and $n$-ary operations. Algebras of $n$-place
functions closed with respect to these compositions were
investigated, for example, in \cite{Sok} and \cite{Yak}.

\section{Preliminaries and notations}

Let $A^n$ be the $n$-th Cartesian product of a set $A$. Any
partial mapping from $A^n$ into $A$ is called a {\it partial
$n$-place function}. The set of all such mappings is denoted by
${\mathcal F}(A^n,A)$. On ${\mathcal F}(A^n,A)$ we define the
superposition (composition) of $n$-place functions
$O:(f,g_1,\ldots,g_n)\mapsto f[g_1\ldots g_n]$ and $n$ binary
compositions $\op{1},\ldots,\op{n}$ putting
\begin{eqnarray}\label{1}
&&f[g_1\ldots g_n](a_1,\ldots,a_n) =
f(g_1(a_1,\ldots,a_n),\ldots,g_n(a_1,\ldots,a_n)),\\[4pt]
\label{2}
&&(f\op{i\,}g)(a_1,\ldots,a_n)=f(a_1,\ldots,a_{i-1},g(a_1,\ldots,a_n),a_{i+1},\ldots,a_n),
\end{eqnarray}
for all $\;f,g,g_1,\ldots,g_n\in {\mathcal F}(A^n,A)$ and
$(a_1,\ldots,a_n)\in A^n$, where left and right side of (\ref{1})
and (\ref{2}) are defined or not defined simultaneously. Since, as
it is not difficult to verify, each composition $\op{i\,}$ is an
associative operation, algebras of the form
$(\Phi;\op{1\,},\ldots,\op{n})$ and
$(\Phi;O,\op{1\,},\ldots,\op{n})$, where
$\Phi\subset\mathcal{F}(A^n,A)$, are called respectively {\it
$(2,n)$-semigroups} and {\it Menger $(2,n)$-semigroups of
$n$-place functions}.

According to the general convention used in the theory of $n$-ary
systems, the sequence $\,x_i,x_{i+1},\ldots,x_j$, where
$i\leqslant j$, \ can be written as $\,x_i^j$ \ (for \ $i>j$ \ it
is the empty symbol). In this convention (\ref{1}) and (\ref{2})
can be written as
\begin{eqnarray*}
&&f[g_1^n](a_1^n) = f(g_1(a_1^n),\ldots,g_n(a_1^n)),\\[4pt]
&&(f\op{i\,}g)(a_1^n)=f(a_1^{i-1},g(a_1^n),a_{i+1}^n).
\end{eqnarray*}

An algebra $(G;o)$ with one $(n+1)$-ary operation $o$ satisfying
the identity
$$
o(o(x_0^n),y_1^n)=o(x_0,o(x_1,y_1^n),\ldots,o(x_n,y_1^n))
$$
is called a {\it Menger algebra of rank $\,n$} (cf. \cite{DT1},
\cite{ST}). Such operation is called {\it superassociative} and by
many authors is written as $o(x_0^n)=x_0[x_1^n]$. In this
convention the above identity has the form
\begin{equation}\label{3}
x_0[x_1^n][y_1^n] = x_0[x_1[y_1^n ]\ldots x_n[y_1^n ]\,] .
\end{equation}
It is clear that a Menger algebra of rank $1$ is an arbitrary
semigroup.

Let $\{\op{1\,},\ldots,\op{n}\}$ be a collection of associative
binary operations defined on $G$. According to \cite{Sok} and
\cite{Yak}, an algebra $(G;\op{1\,},\ldots,\op{n})$ is called a
{\em $(2,n)$-semigroup}. By a {\em Menger $(2,n)$-semigroup} we
mean an algebra $(G;o,\op{1\,},\ldots,\op{n})$, where $(G;o)$ is a
Menger algebra of rank $n$ and $(G;\op{1\,},\ldots,\op{n})$ is a
$(2,n)$-semigroup. Any homomorphism of a (Menger)
$(2,n)$-semigroup onto some (Menger) $(2,n)$-semigroup of
$n$-place functions is called a {\it representation by $n$-place
functions}. A representation is {\it faithful} if it is an
isomorphism.

The symbol $\mu_i(\opp{i_1}{i_s}x_1^s)$, where $x_1,\ldots,x_s\in
G$ and $\op{i_1},\ldots,\op{i_s}$ are binary operations defined on
$G$, denotes an element \
$x_{i_k}\!\opp{i_{k+1}}{i_s}\!x_{k+1}^{s}$ \ if $i=i_k$\ and\
$i\neq i_p$ for all $p<k\leqslant s$. If\ $i\neq i_p$\ for all
$i_p\in\{i_1,\ldots,i_s\}$\ this symbol is empty. For example,
$\mu_1(\op{2}x\op{1\,}y\op{3}z)=y\op{3}z$,
$\mu_2(\op{2}x\op{1\,}y\op{3}z)=x\op{1}y\op{3}z$,
$\mu_3(\op{2}x\op{1\,}y\op{3}z)=z$. The symbol
$\mu_4(\op{2}x\op{1}y\op{3}z)$ is empty.

In \cite{Sok} it is proved that a $(2,n)$-semigroup
$(G;\op{1\,},\ldots,\op{n})$ has a faithful representation by
$n$-place functions if and only if it satisfies the implication
\begin{equation}
\label{4}
\bigwedge\limits_{i=1}^{n}\left(\mu_{i}(\opp{i_{1}}{i_{s}}x_{1}^{s})=
\mu_{i}(\opp{j_{1}}{j_{k}}y_{1}^{k})\right)\longrightarrow
g\opp{i_{1}}{i_{s}}x_{1}^{s}=g\opp{j_{1}}{j_{k}}y_{1}^{k}.
\end{equation}
For Menger $(2,n)$-semigroups the following identities must be
satisfied additionally
\begin{eqnarray}
&&\label{5}(x\op{i\,}y)[z_{1}^{n}]=x[z_{1}^{i-1}y[z_{1}^{n}]\,z_{i+1}^{n}],
\\[4pt]
&&\label{6}x[y_{1}^{n}]\op{i\,}z=x[y_{1}\!\op{i\,}z\ldots
y_{n}\op{i\,}z],
\\[3pt]
&&\label{7}x\opp{i_{1}}{i_{s}}y_{1}^{s}=
x[\mu_{1}(\opp{i_{1}}{i_{s}}y_{1}^{s})\ldots
\mu_{n}(\opp{i_{1}}{i_{s}}y_{1}^{s})],
\end{eqnarray}
where $\{i_{1},\ldots,i_{s}\}=\{1,\ldots,n\}$ and $i=1,\ldots,n$.
In the sequel, any (Menger) $(2,n)$-semigroup satisfying the
condition (\ref{4}) (respectively, (\ref{4}), (\ref{5}), (\ref{6})
and (\ref{7})) will be called  {\it representable}.

Let $\Phi$ be some set of $n$-place functions, i.e.
$\Phi\subset\mathcal{F}(A^{n},A)$. Consider the following three
binary relations on $\Phi$:
\begin{eqnarray*}
&&\chi_{\Phi}=\{(f,g)\in\Phi\times\Phi\,|\,\pr f\subset\pr g\},
\\[4pt]
&&\gamma_{\Phi}=\{(f,g)\in\Phi\times\Phi\,|\,\pr f\cap\pr
g\neq\varnothing\},
\\[4pt]
&&\pi_{\Phi}=\{(f,g)\in\Phi\times\Phi\,|\,\pr f=\pr g\},
\end{eqnarray*}
where $\pr f$ is the domain of $f$, called respectively: {\it
inclusion of domains}, {\it co-definability} and {\it equality of
domains}.

Abstract characterizations of such relations for semigroups of
transformations were studied in \cite{Sch1}, \cite{Sch2},
\cite{Sch3} and for Menger algebras of $n$-place functions in
\cite{Tr1}, \cite{Tr2}, \cite{Tr3}. In this paper these relations
will be characterized in $(2,n)$-semigroups and Menger
$(2,n)$-semigroups of $n$-place functions.

Consider a representable (Menger) $(2,n)$-semigroup
$(G;\op{1\,},\ldots,\op{n})$ (respectively,
$(G;o,\op{1\,},\ldots,\op{n})$) and its representation $P$ by
$n$-place functions. On the set $G$ we define the following three
binary relations:
\begin{eqnarray}
&&\chi_{P}=\{(g_{1},g_{2})\,|\,\pr P(g_{1})\subset\pr P(g_{2})\},  \nonumber\\[4pt]
&&\gamma_{P}=\{(g_{1},g_{2})\,|\,\pr P(g_{1})\cap\pr P(g_{2})\neq\varnothing\},
\nonumber\\[4pt]
&&\pi_{P}=\{(g_{1},g_{2})\,|\,\pr P(g_{1})=\pr P(g_{2})\}.
\nonumber
\end{eqnarray}
It is not difficult to see that $\chi_{P}$ is a quasi-order and
$\pi_{P}$ is an equivalence such that
$\pi_{P}=\chi_{P}\cap\chi_{P}^{-1}$, where
$\chi_{P}^{-1}=\{(b,a)\,|\,(a,b)\in\chi_{P}\}$.

Let $(P_{i})_{i\in I}$ be a family of representations of a
representable $(2,n)$-semi\-group $(G;\op{1\,},\ldots,\op{n})$
(respectively, representable Menger $(2,n)$-semi\-group
$(G;o,\op{1\,},\ldots,\op{n})$) by $n$-place functions defined on
sets $(A_{i})_{i\in I}$ respectively, where the sets $A_{i}$ are
pairwise disjoint. The \textit{sum} of $(P_{i})_{i\in I}$ is the
mapping $P:g\mapsto P(g)$, denoted by $\sum\limits_{i\in I}P_i$,
where $P(g)$ is an $n$-place function on $A=\bigcup\limits_{i\in
I}A_{i}$ such that $P(g)=\bigcup\limits_{i\in I}P_{i}(g)$ for
every $g\in G$. The sum of a family of representations by
$n$-place functions is also a representation by $n$-place
functions and
\begin{equation}\label{8}
  \chi_{P}=\bigcap\limits_{i\in I}\chi_{P_{i}},\quad
  \gamma_{P}=\bigcup\limits_{i\in I}\gamma_{P_{i}},\quad
  \pi_{P}=\bigcap\limits_{i\in I}\pi_{P_{i}}.
\end{equation}

Let $0$ be a zero of a $(2,n)$-semigroup
$(G;\op{1\;},\ldots,\op{n})$ (respectively, Menger
$(2,n)$-semigroup $(G;o,\op{1},\ldots,\op{n})$), i.e.
$0\op{i\,}g=g\op{i\,}0=0$ (respectively, $0\op{i\,}g=g\op{i\,}0=0$
and $\,0[g_{1}^{n}]=g[g_{1}^{i-1}\,0\,g_{i+1}^{n}]=0$) for all
$i=1,\ldots,n$ and $g,g_{1},\ldots,g_{n}\in G$. We say that a
binary relation $\rho\subset G\times G$ is \textit{$0$-reflexive},
if $(g,g)\in\rho$ for all $g\in G\setminus\{0\}$. A symmetric
relation $\rho$ which is reflexive if $0\in\pr\rho$, and
$0$-reflexive if $0\not\in\pr\rho$, is called a
\textit{$0$-quasi-equivalence}.

A binary relation $\Delta$ on a Menger $(2,n)$-semigroup
$(G;o,\op{1},\ldots,\op{n})$ is called:
\begin{itemize}
  \item \textit{$l$-regular}, if
  \begin{eqnarray}
  &&\label{9}x\,\Delta\, y\longrightarrow x[z_{1}^{n}]\;\Delta\; y[z_{1}^{n}], \\[4pt]
  &&\label{10}x\,\Delta\, y\longrightarrow x\op{i\,}z\;\Delta\; y\op{i\,}z
  \end{eqnarray}
  for all $i=1,\ldots,n$ and $x,y,z,z_{1},\ldots,z_{n}\in G$,
  \item \textit{$l$-cancellative}, if
  \begin{eqnarray}
  &&\label{11}x[z_{1}^{n}]\;\Delta\;y[z_{1}^{n}]\longrightarrow x\,\Delta\, y, \\[4pt]
  &&\label{12}x\op{i\,}z\;\Delta\; y\op{i\,}z\longrightarrow x\,\Delta\, y
  \end{eqnarray}
  for all $ i=1,\ldots,n$ and $x,y,z,z_{1},\ldots,z_{n}\in G$,
  \item \textit{$v$-negative}, if
  \begin{eqnarray}
&&\label{13}x[y_{1}^{n}]\;\Delta\;y_{i}, \ i=1,\ldots,n, \\[4pt]
&&\label{14}x\opp{i_{1}}{i_{s}}y_{1}^{s}\;\Delta\;\mu_{j}(\opp{i_{1}}{i_{s}}y_{1}^{s})
  \end{eqnarray}
for all $\,x,y_{1},\ldots,y_{k}\in G$, $\,k=max\{n,s\}$ and
$\,j\in\{i_{1},\ldots,i_{s}\}$.
\end{itemize}
In the case of $(2,n)$-semigroups these relations are defined only
by (\ref{10}), (\ref{12}) and (\ref{14}), respectively.

\section{Projection representable relations on Menger
$(2,n)$-semigroups}

\noindent Let $\mathcal{G}=(G;o,\op{1\,},\ldots,\op{n})$ be a
representable Menger $(2,n)$-semigroup, $\chi$, $\gamma$, $\pi$~--
binary relations on $G$. We say that the triplet
$(\chi,\gamma,\pi)$ is ({\it faithful}) \textit{projection
representable for $\mathcal{G}$}, if there exists such (faithful)
representation $P$ of $\mathcal{G}$ by $n$-place functions for
which $\chi=\chi_{P}$, $\gamma=\gamma_{P}$ and $\pi=\pi_{P}$.
Analogously we define projection representable pairs and separate
relations.

In the sequel, instead of $(g_{1},g_{2})\in\chi$,
$(g_{1},g_{2})\in\gamma$ and $(g_{1},g_{2})\in\pi$ we will write
$\,g_{1}\sqsubset g_{2}$, $\,g_{1}\top g_{2}\,$ and $\,g_{1}\equiv
g_{2}$, respectively.

\begin{Theorem}\label{T1}
A triplet $(\chi,\gamma,\pi)$ of binary relations on $G$ is
projection representable for a representable Menger
$(2,n)$-semigroup $\mathcal{G}$ if and only if the following
conditions are satisfied:

 $(a)$ \ $\chi$ is an
$l$-regular and $v$-negative quasi-order,

$(b)$ \ $\gamma$ is an $l$-cancellative $0$-quasi-equivalence,

$(c)$ \ $\pi=\chi\cap\chi^{-1}$ \ and
\begin{equation}\label{15}
h_{1}\top h_{2}\,\wedge\,h_{1}\sqsubset g_{1}\,\wedge\,
h_{2}\sqsubset g_{2}\longrightarrow g_{1}\top g_{2}\,
\end{equation}

\hspace*{4mm} for all $\,h_{1},h_{2},g_{1},g_{2}\in G$.
\end{Theorem}
\begin{proof}{\it Necessity.} Let
$(\Phi;O,\op{1\;},\ldots,\op{n})$ be a Menger $(2,n)$-semigroup of
$n$-place functions determined on the set $A$. Let us show that
the triplet $(\chi_{\Phi},\gamma_{\Phi},\pi_{\Phi})$ satisfies all
the conditions of the theorem.

At first we prove the condition $(a)$. The relation $\chi_{\Phi}$
is obviously a quasi-order. Let $f,g,h_{1},\ldots,h_{n}\in\Phi$
and $(f,g)\in\chi_{\Phi}$, i.e. $\pr f\subset\pr g$. Suppose that
$\bar{a}\in\pr f[h_{1}^{n}]$ for some $\bar{a}\in A^{n}$. Then
$\{f[h_{1}^{n}](\bar{a})\}\neq\varnothing$, i.e.
$\{f(h_{1}(\bar{a}),\ldots h_{n}(\bar{a}))\}\neq\varnothing$. Thus
$(h_{1}(\bar{a}),\ldots h_{n}(\bar{a}))\in\pr f$ and, in the
consequence, $(h_{1}(\bar{a}),\ldots h_{n}(\bar{a}))\in\pr g$.
Therefore $\{g(h_{1}(\bar{a}),\ldots
h_{n}(\bar{a}))\}\neq\varnothing$, whence
$\{g[h_{1}^{n}](\bar{a})\}\neq\varnothing$, i.e. $\bar{a}\in\pr
g[h_{1}^{n}]$. So, $\pr f[h_{1}^{n}]\subset\pr g[h_{1}^{n}]$,
which implies $(f[h_{1}^{n}],g[h_{1}^{n}])\in\chi_{\Phi}$.
Similarly we can prove that for all $f,g,h\in\Phi$ and
$i=1,\ldots,n$, from $(f,g)\in\chi_{\Phi}$ it follows
$(f\op{i\,}h\, ,g\op{i\,}h)\in\chi_{\Phi}$. This means that the
relation $\chi_{\Phi}$ is $l$-regular. The proof of the
$v$-negativity is analogous.

To prove $(b)$ let $\Theta$ be a zero of a Menger
$(2,n)$-semigroup $(\Phi;O,\op{1\;},\ldots,\op{n})$. If
$\Theta\neq\varnothing$, then $\pr\Theta\neq\varnothing$, whence
$(\Theta,\Theta)\in\gamma_{\Phi}$. Thus
$\Theta\in\pr\gamma_{\Phi}$. So, in this case $\gamma_{\Phi}$ is
reflexive. For $\Theta=\varnothing$ we have
$\pr\Theta=\varnothing$. Therefore
$\Theta\not\in\pr\gamma_{\Phi}$, i.e. $(f,f)\in\gamma_{\Phi}$ for
every $f\neq\Theta$. Hence $\gamma_{\Phi}$ is $\Theta$-reflexive.
Since $\gamma_{\Phi}$ is symmetric, the above means that
$\gamma_{\Phi}$ is a $\Theta$-quasi-equivalence.

Suppose now that $(f[h_{1}^n],g[h_{1}^n])\in\gamma_{\Phi}$ for
some $f,g\in\Phi$, $h_{1}^{n}\in\Phi^{n}$. Then $\pr
f[h_{1}^n]\cap\pr g[h_{1}^n]\neq\varnothing$, i.e. there exists
$\bar{a}\in A^{n}$ such that $\bar{a}\in\pr f[h_{1}^{n}]$ and
$\bar{a}\in\pr g[h_{1}^{n}]$. Therefore
$\{f[h_{1}^{n}](\bar{a})\}\neq\varnothing$ and
$\{g[h_{1}^{n}](\bar{a})\}\neq\varnothing$. Thus
$\{f(h_{1}(\bar{a}),\ldots h_{n}(\bar{a}))\}\neq\varnothing$ and
$\{g(h_{1}(\bar{a}),\ldots h_{n}(\bar{a}))\}\neq\varnothing$,
which shows that $(h_{1}(\bar{a}),\ldots,h_{n}(\bar{a}))\in\pr
f\cap\pr g$. So, $(f,g)\in\gamma_{\Phi}$. Analogously, for
$f,g,h\in\Phi$, $i=1,\ldots,n$, from
$(f\op{i\,}h,g\op{i\,}h)\in\gamma_{\Phi}$ it follows
$(f,g)\in\gamma_{\Phi}$. So, $\gamma_{\Phi}$ is $l$-cancellative.

Since in $(c)$ the first condition is obvious, we prove
$(\ref{15})$ only. For this let $(h_{1},h_{2})\in\gamma_{\Phi}$,
$(h_{1},g_{1})\in\chi_{\Phi}$ and $(h_{2},g_{2})\in\chi_{\Phi}$
for some $h_{1},h_{2},g_{1},g_{2}\in\Phi$. Then $\,\pr
h_{1}\cap\pr h_{2}\neq\varnothing$, $\,\pr h_{1}\subset\pr g_{1}$
and $\,\pr h_{2}\subset\pr g_{2}$, whence $\varnothing\neq\pr
h_{1}\cap\pr h_{2}\subset\pr g_{1}\cap\pr g_{2}$. Thus $\,\pr
g_{1}\cap\pr g_{2}\neq\varnothing$, i.e.
$(g_{1},g_{2})\in\gamma_{\Phi}$, which proves (\ref{15}) and
completes the proof of the necessity of the conditions formulated
in the theorem.
 \end{proof}

To prove the sufficiency of these conditions we must introduce
some additional constructions. Consider the triplet
$(\chi,\gamma,\pi)$ of binary relations on a representable Menger
$(2,n)$-semigroup $\mathcal{G}=(G;o,\op{1\;},\ldots,\op{n})$
satisfying all the conditions of the theorem. Let
$e_{1},\ldots,e_{n}$ be pairwise different elements not belonging
to $G$. For all $x_{1},\ldots,x_{s}\in G$, $i=1,\ldots,n$, and
operations $\,\op{i_{1}},\ldots,\op{i_{s}}$ defined on $G$ by
$\mu_{i}^{*}(\opp{i_{1}}{i_{s}}x_{1}^{s})$ we denote an element of
$G^{\,\ast}=G\cup\{e_{1},\ldots,e_{n}\}$ such that
$$
\mu_{i}^{*}(\opp{i_{1}}{i_{s}}x_{1}^{s})=\left\{
\begin{array}{cl}
\mu_{i}(\opp{i_{1}}{i_{s}}x_{1}^{s}),&\mbox{if }\ i\in\{i_{1},\ldots,i_{s}\},\\[4mm]
  e_{i}, & \mbox{if }\ i\not\in\{i_{1},\ldots,i_{s}\}.
\end{array}
\right.
$$

Consider the set
$\frak{A}^*=G^n\cup\frak{A}_{0}\cup\{(e_{1},\ldots,e_{n})\}$,
where $\frak{A}_{0}$\label{Tr**} is the collection of all
$n$-tuples $(x_1,\ldots,x_n)\in (G^*)^n$ for which there exists
$y_1,\ldots,y_s\in G$ and $\,i_1,\ldots,i_n\in\{1,\ldots,n\}$ such
that $x_i=\mu_i^*(\opp{i_1}{i_s}y_1^s)$. Let $(h_1,h_2)\in G^2$ be
fixed. For each $g\in G$ we define a partial $n$-place function
$P_{(h_{1},\,h_{2})}(g):\frak{A}^*\rightarrow G$ such that
\[
x_{1}^{n}\in\pr P_{(h_{1},\,h_{2})}(g)\longleftrightarrow
 \left\{
\begin{array}{ll}
h_{1}\sqsubset g[x_{1}^{n}]\,\vee\,h_{2}\sqsubset g[x_{1}^{n}]
 & \mbox{if }\ x_{1}^{n}\in G^{n}, \\[4pt]
h_{1}\sqsubset g\,\vee\, h_{2}\sqsubset g
 & \mbox{if }\ x_{1}^{n}=e_{1}^{n}, \\[4pt]
h_{1}\sqsubset g\opp{i_{1}}{i_{s}}y_{1}^{s}\,\vee\,h_{2}\sqsubset
g\opp{i_{1}}{i_{s}}y_{1}^{s} & \mbox{if }\
  x_{i}=\mu_{i}^{*}(\opp{i_{1}}{i_{s}}y_{1}^{s}),\\
  & i=1,\ldots,n,\mbox{ for}  \\
  & \mbox{some }y_{1}^{s}\in G^{s}\; \mbox{and }\\
  & i_{1}\ldots,i_{s}\in\{1,\ldots,n\}.
\end{array}
  \right.
\]
For $\,x_{1}^{n}\in\pr P_{(h_{1},\,h_{2})}(g)$ we put
\begin{equation}\label{16}
P_{(h_{1},\,h_{2})}(g)(x_{1}^{n})= \left\{
\begin{array}{ll}
  g[x_{1}^{n}] & \mbox{if }\ x_{1}^{n}\in G^{n}, \\[4pt]
  g & \mbox{if }\ x_{1}^{n}=e_{1}^{n}, \\[4pt]
  g\opp{i_{1}}{i_{s}}y_{1}^{s} & \mbox{if }\
  x_{i}=\mu_{i}^{*}(\opp{i_{1}}{i_{s}}y_{1}^{s}),\\
  & i=1,\ldots,n,\mbox{ for}  \\
  & \mbox{some }y_{1}^{s}\in G^{s}\  \mbox{and }\\
  & i_{1}\ldots,i_{s}\in\{1,\ldots,n\}.
\end{array}
  \right.
\end{equation}
Let us show that $P_{(h_{1},\,h_{2})}$ is a representation of
$\mathcal{G}$ by $n$-place functions.

\begin{Proposition}\label{P1}
The function $P_{(h_{1},\,h_{2})}(g)$ is single-valued.
\end{Proposition}
\begin{proof}
Let $x_1^n\in\pr P_{(h_{1},\,h_{2})}(g)$, where $g,h_1,h_2\in G$
are fixed. Since for $x_{1}^{n}\in G^{n}$ and
$x_{1}^{n}=e_{1}^{n}$ the value of
$P_{(h_{1},\,h_{2})}(g)(x_{1}^{n})$ is uniquely determined, we
verify only the case when
$\,x_{i}=\mu_{i}^{*}(\opp{i_{1}}{i_{s}}y_{1}^{s})$,
$i=1,\ldots,n$, for some $y_{1}^{s}\in G^{s}$. If for some
$z_{1}^{k}\in G^{k}$ and $j_{1},\ldots,j_{k}\in\{1,\ldots,n\}$ we
have also $x_{i}=\mu_{i}^{*}(\opp{j_{1}}{j_{k}}z_{1}^{k})$,
$i=1,\ldots,n$, then $\mu_{i}(\opp{i_{1}}{i_{s}}y_{1}^{s})=
\mu_{i}(\opp{j_{1}}{j_{k}}z_{1}^{k})$ for every $i=1,\ldots,n$,
which, according to (\ref{4}), implies
$g\opp{i_{1}}{i_{s}}y_{1}^{s}=g\opp{j_{1}}{j_{k}}z_{1}^{k}$. This
means that also in this case $P_{(h_{1},\,h_{2})}(g)(x_{1}^{n})$
is uniquely determined. Thus, the function
$P_{(h_{1},\,h_{2})}(g)$ is single-valued.
\end{proof}

\begin{Proposition}\label{P2}
For all $g,g_{1},\ldots,g_{n},h_{1},h_{2}\in G$ we have
\[
P_{(h_{1},\,h_{2})}(g[g_{1}^{n}])=
P_{(h_{1},\,h_{2})}(g)[P_{(h_{1},\,h_{2})}(g_{1})\ldots
P_{(h_{1},\,h_{2})}(g_{n})].
\]
\end{Proposition}
\begin{proof} Let $g,g_{1},\ldots,g_{n}\in G$ and
$x_{1}^{n}\in\pr P_{(h_{1},\,h_{2})}(g[g_{1}^{n}])$. If
$x_{1}^{n}\in G^{n}$, then
\[
h_{1}\sqsubset g[g_{1}^{n}][x_{1}^{n}]\,\vee\, h_{2}\sqsubset
  g[g_{1}^{n}][x_{1}^{n}],
\]
whence, applying the superassociativity (\ref{3}), we obtain
\begin{equation}\label{17}
h_{1}\sqsubset g[g_{1}[x_{1}^{n}]\ldots g_{n}[x_{1}^{n}]]\,\vee\,
h_{2}\sqsubset g[g_{1}[x_{1}^{n}]\ldots g_{n}[x_{1}^{n}]].
\end{equation}
This together with the $v$-negativity of $\chi$ implies
\begin{equation}\label{18}
  h_{1}\sqsubset g_{i}[x_{1}^{n}]\,\vee\, h_{2}\sqsubset
  g_{i}[x_{1}^{n}],\ i=1,\ldots,n.
\end{equation}
From (\ref{17}) it follows that
$(g_{1}[x_{1}^{n}],\ldots,g_{n}[x_{1}^{n}])\in\pr
P_{(h_{1},\,h_{2})}(g)$, from (\ref{18}) that $x_{1}^{n}\in
P_{(h_{1},\,h_{2})}(g_{i})$, $i=1,\ldots,n$. So, if $\,x_1^n\in
G^n$, then
\begin{equation}\label{19}
  x_1^n\in\pr
  P_{(h_{1},\,h_{2})}(g[g_1^n])\longleftrightarrow\left\{
\begin{array}{l}
  (g_{1}[x_{1}^{n}],\ldots,g_{n}[x_{1}^{n}])\in\pr
P_{(h_{1},\,h_{2})}(g), \\[4pt]
  \bigwedge\limits_{i=1}^{n}x_{1}^{n}\in
  P_{(h_{1},\,h_{2})}(g_{i}).
\end{array}
  \right.
\end{equation}
Analogously we can verify that
\begin{equation}\label{20}
  e_1^n\in\pr
  P_{(h_{1},\,h_{2})}(g[g_1^n])\longleftrightarrow\left\{
\begin{array}{l}
  (g_{1},\ldots,g_{n})\in\pr
P_{(h_{1},\,h_{2})}(g), \\[4pt]
  \bigwedge\limits_{i=1}^{n}e_{1}^{n}\in
  P_{(h_{1},\,h_{2})}(g_{i}).
\end{array}
  \right.
\end{equation}

Now let $x_{i}=\mu_{i}^{*}(\opp{i_{1}}{i_{s}}y_{1}^{s})$,
$i=1,\ldots,n$, for some $i_1,\ldots,i_s\in\{1,\ldots,n\}$ and
$y_{1}^{s}\in G^{s}$. Then $x_1^n\in\pr
  P_{(h_{1},\,h_{2})}(g[g_1^n])$ implies
\[
h_1\sqsubset g[g_1^n]\opp{i_1}{i_s}y_1^s\;\vee\; h_2\sqsubset
  g[g_1^n]\opp{i_1}{i_s}y_1^s,
\]
which, by (\ref{6}), is equivalent to
\begin{equation}\label{21}
  h_1\sqsubset g[g_1\!\opp{i_1}{i_s}y_1^s\ldots
  g_n\!\opp{i_1}{i_s}y_1^s]\;\vee\;h_2\sqsubset g[g_1\!\opp{i_1}{i_s}y_1^s\ldots
  g_n\!\opp{i_1}{i_s}y_1^s].
\end{equation}
From this, applying the $v$-negativity of $\chi$, we
obtain
\begin{equation}\label{22}
  h_1\sqsubset g_i\opp{i_1}{i_s}y_1^s\;\vee\; h_2\sqsubset
  g_i\opp{i_1}{i_s}y_1^s
\end{equation}
for every $\,i=1,\ldots,n$.

The condition (\ref{21}) is equivalent to
$(g_1\opp{i_1}{i_s}y_1^s,\ldots,g_n\opp{i_1}{i_s}y_1^s)\in\pr
P_{(h_{1},\,h_{2})}(g)$. The condition (\ref{22}) shows that
$x_1^n\in\pr P_{(h_{1},\,h_{2})}(g_{i})$ for every $i=1,\ldots,n$,
where $x_i=\mu_i^*(\opp{i_1}{i_s}y_1^s)$, $\,i=1,\ldots,n$. So,
\begin{equation}\label{23}
  x_1^n\in\pr
  P_{(h_{1},\,h_{2})}(g[g_1^n])\longleftrightarrow\left\{
\begin{array}{l}
(g_1\opp{i_1}{i_s}y_1^s,\ldots,g_n\opp{i_1}{i_s}y_1^s)\in\pr P_{(h_{1},\,h_{2})}(g),\\[4mm]
  \bigwedge\limits_{i=1}^{n}x_1^n\in\pr P_{(h_{1},\,h_{2})}(g_i),
\end{array}
  \right.
\end{equation}
where $\,x_i=\mu_i^*(\opp{i_1}{i_s}y_1^s)$, $\,i=1,\ldots,n$.

Let $\,x_1^n\in\pr P_{(h_{1},\,h_{2})}(g[g_1^n])$. If $\,x_1^n\in
G^n$, then, according to (\ref{16}) and (\ref{19}), we have
\[\arraycolsep=.5mm
\begin{array}{rl}
P_{(h_{1},\,h_{2})}(g[g_1^n])(x_1^n)&=g[g_1^n][x_1^n]=
g[g_1[x_1^n]\ldots g_n[x_1^n]] \\[4pt]
&=P_{(h_{1},\,h_{2})}(g)(g_1[x_1^n],\ldots, g_n[x_1^n]) \\[4pt]
&=P_{(h_{1},\,h_{2})}(g)\left(P_{(h_{1},\,h_{2})}(g_1)(x_1^n),\ldots,
P_{(h_{1},\,h_{2})}(g_n)(x_1^n)\right) \\[4pt]
&=P_{(h_{1},\,h_{2})}(g)\left[P_{(h_{1},\,h_{2})}(g_1)\ldots
P_{(h_{1},\,h_{2})}(g_n)\right](x_1^n).
\end{array}
\]

Similarly, we can prove that
\[
P_{(h_{1},\,h_{2})}(g[g_1^n])(e_1^n)=\left[P_{(h_{1},\,h_{2})}(g_1)\ldots
P_{(h_{1},\,h_{2})}(g_n)\right](e_1^n)
\]
for $e_1^n\in\pr P_{(h_{1},\,h_{2})}(g[g_1^n])$.

If $x_1^n\in\pr P_{(h_{1},\,h_{2})}(g[g_1^n])$, where
$x_{i}=\mu_{i}^{*}(\opp{i_{1}}{i_{s}}y_{1}^{s})$, $i=1,\ldots,n$,
for some $y_{1}^{s}\in G^{s}$,
$\,i_1,\ldots,i_s\in\{1,\ldots,n\}$, then, according to (\ref{16})
and (\ref{23}), we obtain
\[\arraycolsep=.5mm
\begin{array}{rl}
P_{(h_{1},\,h_{2})}(g[g_1^n])(x_1^n)&=g[g_1^n]\opp{i_{1}}{i_{s}}y_{1}^{s}
=g[g_1\opp{i_{1}}{i_{s}}y_{1}^{s}\ldots
g_n\opp{i_{1}}{i_{s}}y_{1}^{s}]\\
&=P_{(h_{1},\,h_{2})}(g)(g_1\opp{i_{1}}{i_{s}}y_{1}^{s},\ldots,
g_n\opp{i_{1}}{i_{s}}y_{1}^{s}) \\[5pt]
&=P_{(h_{1},\,h_{2})}(g)\left(P_{(h_{1},\,h_{2})}(g_1)(x_1^n),\ldots,
P_{(h_{1},\,h_{2})}(g_n)(x_1^n)\right) \\[4pt]
&=P_{(h_{1},\,h_{2})}(g)\left[P_{(h_{1},\,h_{2})}(g_1)\ldots
P_{(h_{1},\,h_{2})}(g_n) \right](x_1^n).
\end{array}
\]
The proof is complete.
\end{proof}

\begin{Proposition}\label{P3}
For all $g_1,g_2,h_1,h_2\in G$ and $\,i=1,\ldots,n$ we have
\[
P_{(h_{1},\,h_{2})}(g_1\op{i\,}g_2)=
P_{(h_{1},\,h_{2})}(g_1)\op{i\,}P_{(h_{1},\,h_{2})}(g_2).
\]
\end{Proposition}
\begin{proof} Let $x_1^n\in\pr P_{(h_{1},\,h_{2})}(g_1\op{i\,}g_2)$. If
$x_1^n\in G^n$, then
\[
  h_1\sqsubset(g_1\op{i\,}g_2)[x_1^n]\,\vee\,
  h_2\sqsubset(g_1\op{i\,}g_2)[x_1^n],
\]
which, by (\ref{5}), is equivalent to
\begin{equation}\label{24}
  h_1\sqsubset g_1[x_{1}^{i-1}g_2[x_1^n]x_{i+1}^{n}]\,\vee\, h_2\sqsubset
  g_1[x_{1}^{i-1}g_2[x_1^n]x_{i+1}^{n}].
\end{equation}
This, according to the $v$-negativity of $\chi$, implies
\begin{equation}\label{25}
h_1\sqsubset g_2[x_1^n]\,\vee\, h_2\sqsubset g_2[x_1^n].
\end{equation}
The condition (\ref{24}) means that
$(x_{1}^{i-1},g_2[x_1^n],x_{i+1}^{n})\in\pr
P_{(h_{1},\,h_{2})}(g_1)$. From (\ref{25}) we obtain $x_1^n\in\pr
P_{(h_{1},\,h_{2})}(g_2)$. So, for $x_1^n\in G^n$ we have
\begin{equation}\label{26} x_1^n\in\pr
P_{(h_{1},\,h_{2})}(g_1\op{i\,}g_2)\longleftrightarrow\left\{
\begin{array}{l}
  (x_{1}^{i-1},g_2[x_1^n],x_{i+1}^{n})\in\pr
P_{(h_{1},\,h_{2})}(g_1) \\[4pt]
  x_1^n\in\pr
P_{(h_{1},\,h_{2})}(g_2).
\end{array}
\right.
\end{equation}

Consider now the case when $x_1^n=e_1^n$. In this case
$e_1^n\in\pr P_{(h_{1},\,h_{2})}(g_1\op{i\,}g_2)$ means, by
(\ref{17}), that
\begin{equation}\label{27}
h_1\sqsubset g_1\op{i\;}g_2\,\vee\, h_2\sqsubset g_1\op{i\;}g_2.
\end{equation}
Because $\,g_1\op{i\;}g_2\sqsubset\mu_i(\op{i\,}g_2)=g_2\,$, by
the $v$-negativity of $\chi$, the above condition gives
\begin{equation}\label{28}
  h_1\sqsubset g_2\,\vee\, h_2\sqsubset g_2.
\end{equation}
But $\mu_{i}^{*}(\op{i\,}g_2)=\mu_{i}(\op{i\,}g_2)=g_2$ and
$\mu_{k}^{*}(\op{i\,}g_2)=e_k$ for
$k\in\{1,\ldots,n\}\setminus\{i\}$, so, (\ref{27}) implies
$(e_{1}^{i-1},g_2,e_{i+1}^{n})\in\pr P_{(h_{1},\,h_{2})}(g_1)$. On
the other hand, from (\ref{28}) it follows $\,e_1^n\in\pr
P_{(h_{1},\,h_{2})}(g_2)$. Therefore
\begin{equation}\label{29}
e_1^n\in\pr
P_{(h_{1},\,h_{2})}(g_1\op{i\,}g_2)\longleftrightarrow\left\{
\begin{array}{l}
  (e_{1}^{i-1},g_2,e_{i+1}^{n})\in\pr P_{(h_{1},\,h_{2})}(g_1) \\[4pt]
  e_1^n\in\pr P_{(h_{1},\,h_{2})}(g_2).
\end{array}
\right.
\end{equation}

In the third case when
$x_{i}=\mu_{i}^{*}(\opp{i_{1}}{i_{s}}y_{1}^{s})$,
$\,i=1,\ldots,n$, for some $y_{1}^{s}\in G^{s}$,
$i_1,\ldots,i_s\in\{1,\ldots,n\}$, from $x_1^n\in\pr
P_{(h_{1},\,h_{2})}(g_1\op{i}g_2)$ we conclude
\begin{equation}\label{30}
  h_1\sqsubset(g_1\op{i\,}g_2)\opp{i_1}{i_s}y_1^s\;\vee\;
  h_2\sqsubset(g_1\op{i\,}g_2)\opp{i_1}{i_s}y_1^s.
\end{equation}
Since $\chi$ is $v$-negative, we have
$(g_1\op{i\;}g_2)\opp{i_1}{i_s}y_1^s\sqsubset\mu_i(\op{i\;}g_2\opp{i_1}{i_s}y_1^s)
=g_2\opp{i_1}{i_s}y_1^s$, which means that (\ref{30}) can be
written in the form
\begin{equation}\label{31}
  h_1\sqsubset g_2\opp{i_1}{i_s}y_1^s\;\vee\; h_2\sqsubset
  g_2\opp{i_1}{i_s}y_1^s.
\end{equation}
But $\mu_i^*(\op{i}g_2\opp{i_1}{i_s}y_1^s)=
\mu_i(\op{i\,}g_2\opp{i_1}{i_s}y_1^s)=g_2\opp{i_1}{i_s}y_1^s$ and
$\mu_k^*(\op{i\,}g_2\opp{i_1}{i_s}y_1^s)=\mu_k^*(\opp{i_1}{i_s}y_1^s)$
for $k\in\{1,\ldots,n\}\setminus\{i\}$. This, together with the
condition (\ref{30}), proves
$(x_{1}^{i-1},g_2\opp{i_1}{i_s}y_1^s,x_{i+1}^{n})\in\pr
P_{(h_{1},\,h_{2})}(g_1).$ Similarly, from (\ref{31}) we can
deduce $x_1^n\in\pr P_{(h_{1},\,h_{2})}(g_2).$ Therefore
\[
x_1^n\in\pr
P_{(h_{1},\,h_{2})}(g_1\op{i\;}g_2)\longleftrightarrow\left\{
\begin{array}{l}
  (x_{1}^{i-1},g_2\opp{i_1}{i_s}y_1^s,x_{i+1}^{n})\in\pr
P_{(h_{1},\,h_{2})}(g_1) \\[4pt]
  x_1^n\in\pr
P_{(h_{1},\,h_{2})}(g_2),
\end{array}
\right.
\]
where $x_i=\mu_i^*(\opp{i_1}{i_s}y_1^s)$, $\,i=1,\ldots,n$.

Let $x_1^n\in\pr P_{(h_{1},\,h_{2})}(g_1\op{i}g_2)$. If $x_1^n\in
G^n$, then, according to (\ref{16}) and (\ref{26}), we have
\[\arraycolsep=.5mm
\begin{array}{rl}
P_{(h_{1},\,h_{2})}(g_1\op{i\;}g_2)(x_1^n)&=(g_1\op{i\;}g_2)[x_1^n]=
g_1[x_{1}^{i-1}g_2[x_1^n]x_{i+1}^{n}] \\[5pt]
&=P_{(h_{1},\,h_{2})}(g_1)\left(x_{1}^{i-1},g_2[x_1^n],x_{i+1}^{n}\right)
 \\[5pt]
&=P_{(h_{1},\,h_{2})}(g_1)\left(x_{1}^{i-1},P_{(h_{1},\,h_{2})}(g_2)(x_1^n),x_{i+1}^{n}\right)
 \\[5pt]
&=P_{(h_{1},\,h_{2})}(g_1)\op{i\,}P_{(h_{1},\,h_{2})}(g_2)(x_1^n).
\end{array}
\]
If $\,x_1^n=e_1^n,\,$ then, analogously as in the previous case,
using (\ref{16}) and (\ref{29}) we
obtain
\[
P_{(h_{1},\,h_{2})}(g_1\op{i\,}g_2)(e_1^n)=
P_{(h_{1},\,h_{2})}(g_1)\op{i\,}P_{(h_{1},\,h_{2})}(g_2)(e_1^n).
\]

Similarly, in the case when
$x_{i}=\mu_{i}^{*}(\opp{i_{1}}{i_{s}}y_{1}^{s})$,
$\,i=1,\ldots,n$, for some $y_{1}^{s}\in G^{s}$,
$i_1,\ldots,i_s\in\{1,\ldots,n\}$, we have
\[\arraycolsep=.5mm
\begin{array}{rl}
P_{(h_{1},\,h_{2})}(g_1\op{i\;}g_2)(x_1^n)&
=(g_1\op{i\,}g_2)\opp{i_{1}}{i_{s}}y_{1}^{s} \\[4pt]
&=P_{(h_{1},\,h_{2})}(g_1)(x_{1}^{i-1},g_2\opp{i_{1}}{i_{s}}y_{1}^{s},x_{i+1}^{n})
 \\[5pt]
&=P_{(h_{1},\,h_{2})}(g_1)\left(x_{1}^{i-1},P_{(h_{1},\,h_{2})}(g_2)(x_1^n),x_{i+1}^{n}\right)
 \\[5pt]
&=P_{(h_{1},\,h_{2})}(g_1)\op{i}P_{(h_{1},\,h_{2})}(g_2)(x_1^n).
\end{array}
\]
This completes our proof.
\end{proof}

Basing on these propositions we are able to prove the sufficiency
of the conditions of Theorem~\ref{T1}.

\medskip
{\it Sufficiency}. Let the triplet $(\chi,\gamma,\pi)$ of binary
relations on a representable Menger $(2,n)$-semigroup
$\mathcal{G}=(G;o,\op{1\,},\ldots,\op{n})$ satisfies all the
conditions of the theorem. Then, as it follows from
Propositions~\ref{P1}--\ref{P3}, for all $h_1,h_2\in G$, the
mapping $P_{(h_{1},\,h_{2})}$ is a representation of $\mathcal{G}$
by $n$-place functions. Consider the family of representations
$P_{(h_{1},\,h_{2})}$ such that ${(h_{1},\,h_{2})\in\gamma}$. Let
$P$ be the sum of this family, i.e.
$P=\!\sum\limits_{(h_{1},\,h_{2})\in\gamma}\!\!P_{(h_{1},\,h_{2})}$.
Of course, $P$ is a representation of $\mathcal{G}$ by $n$-place
functions. Let us show that $\chi=\chi_P$, $\gamma=\gamma_P$ and
$\pi=\pi_P$.

Let $(g_1,g_2)\in\chi_P$. Then, according to (\ref{8}), we have
$(g_1,g_2)\in\chi_{(h_{1},\,h_{2})}$\footnote{\,$\chi_{(h_{1},\,h_{2})}$
denotes this quasi-order which corresponds to the representation
$P_{(h_{1},\,h_{2})}$. Analogously are defined
$\gamma_{(h_{1},\,h_{2})}$ and $\pi_{(h_{1},\,h_{2})}$.} for all
$(h_1,h_2)\in\gamma$, \ i.e.
\[
(\forall(h_1,h_2)\in\gamma)\,\left(\pr
P_{(h_{1},\,h_{2})}(g_1)\subset\pr P_{(h_{1},\,h_{2})}(g_2)
\right),
\]
which is equivalent to
\[
(\forall(h_1,h_2)\in\gamma)(\forall x_1^n)\,\left(x_1^n\in\pr
P_{(h_{1},\,h_{2})}(g_1)\longrightarrow x_1^n\in\pr
P_{(h_{1},\,h_{2})}(g_2)\right).
\]
From this, for $x_1^n=e_1^n$, we obtain
\[
(\forall(h_1,h_2)\in\gamma)\,\left(e_1^n\in\pr
P_{(h_{1},\,h_{2})}(g_1)\longrightarrow e_1^n\in\pr
P_{(h_{1},\,h_{2})}(g_2)\right),
\]
which means that
\[
(\forall(h_1,h_2)\in\gamma)\;\left(h_1\sqsubset g_1\,\vee\,
h_2\sqsubset g_1\longrightarrow h_1\sqsubset g_2\,\vee\,
h_2\sqsubset g_2\right).
\]
Let $g_1\neq 0$. Then $g_1\top\,g_1$ and the above implication
gives $g_1\sqsubset g_1\longrightarrow g_1\sqsubset g_2$. This
proves $(g_1,g_2)\in\chi$ because $\chi$ is reflexive. If $g_1=0$,
then $0=0[g_2\ldots g_2]\sqsubset g_2$, by the $v$-negativity of
$\chi$. Hence $(0,g_2)\in\chi$. So, $(g_1,g_2)\in\chi$, i.e.
$\chi_P\subset\chi$.

Conversely, let $(g_1,g_2)\in\chi$, $(h_1,h_2)\in\gamma$ and
$x_1^n\in\pr P_{(h_{1},\,h_{2})}(g_1)$. If $x_1^n\in G^n$, then
$h_1\sqsubset g_1[x_1^n]\,\vee\,h_2\sqsubset g_1[x_1^n]$. Since
the $l$-regularity of $\chi$ together with $g_1\sqsubset g_2$
implies $g_1[x_1^n]\sqsubset g_2[x_1^n]$, from the above we
conclude $h_1\sqsubset g_2[x_1^n]\vee h_2\sqsubset g_2[x_1^n]$,
i.e. $x_1^n\in\pr P_{(h_{1},\,h_{2})}(g_2)$. Similarly, in the
case $x_1^n=e_1^n$, from $\,e_1^n\in\pr
P_{(h_{1},\,h_{2})}(g_1)\,$ it follows $\,e_1^n\in\pr
P_{(h_{1},\,h_{2})}(g_2)$. In the case when
$x_{i}=\mu_{i}^{*}(\opp{i_{1}}{i_{s}}y_{1}^{s})$,
$\,i=1,\ldots,n$, for some $y_{1}^{s}\in G^{s}$,
$i_1,\ldots,i_s\in\{1,\ldots,n\}$, applying the $l$-regularity of
$\chi$ to $g_1\sqsubset g_2$, we obtain
$g_1\opp{i_{1}}{i_{s}}y_{1}^{s}\sqsubset
g_2\opp{i_{1}}{i_{s}}y_{1}^{s}$, whence, in view of $h_1\sqsubset
g_1\opp{i_{1}}{i_{s}}y_{1}^{s}\,\vee\,h_2\sqsubset
g_1\opp{i_{1}}{i_{s}}y_{1}^{s}$, we obtain $h_1\sqsubset
g_2\opp{i_{1}}{i_{s}}y_{1}^{s}\,\vee\,h_2\sqsubset
g_2\opp{i_{1}}{i_{s}}y_{1}^{s}$. Therefore $x_1^n\in\pr
P_{(h_{1},\,h_{2})}(g_2)$, which proves $\pr
P_{(h_{1},\,h_{2})}(g_1)\subset\pr P_{(h_{1},\,h_{2})}(g_2)$ for
all $(h_{1},\,h_{2})\in\gamma$. Thus $(g_1,g_2)\in\chi_P$, i.e.
$\chi\subset\chi_P$. Consequently, $\chi=\chi_P$. This, together
with the condition $(c)$ formulated in the theorem, gives
$\pi=\chi\cap\chi^{-1}=\chi_{P}\cap\chi_{P}^{-1}=\pi_P$. So,
$\pi=\pi_P$.

Now let $(g_1,g_2)\in\gamma_P$. Then, according to (\ref{8}), we
have $(g_1,g_2)\in\gamma_{(h_{1},\,h_{2})}$ for some
$(h_1,h_2)\in\gamma$, \ i.e.
\[
(\exists(h_1,h_2)\in\gamma)\,\left(\pr
P_{(h_{1},\,h_{2})}(g_1)\cap\pr
P_{(h_{1},\,h_{2})}(g_2)\neq\varnothing\right),
 \]
which is equivalent to
\[
(\exists(h_1,h_2)\in\gamma)(\exists x_1^n)\,\left(x_1^n\in\pr
P_{(h_{1},\,h_{2})}(g_1)\,\wedge\, x_1^n\in\pr
P_{(h_{1},\,h_{2})}(g_2)\right).
\]
This, for $x_1^n\in G^n$ implies $h_1\sqsubset g_1[x_1^n]\,\vee\,
h_2\sqsubset g_1[x_1^n]$ and $h_1\sqsubset g_2[x_1^n]\,\vee\,
h_2\sqsubset g_2[x_1^n]$. From the above, in view of
$h_1\top\,h_2$ and (\ref{15}), we obtain
$g_1[x_1^n]\,\top\,g_2[x_1^n]$, whence, applying the
$l$-cancellativity of $\gamma$, we get $g_1\top\,g_2$, i.e.
$(g_1,g_2)\in\gamma$.

In the similar way, we can see that  in the case $x_1^n=e_1^n$ the
condition $(g_1,g_2)\in\gamma$ also holds.

If $\,x_{i}=\mu_{i}^{*}(\opp{i_{1}}{i_{s}}y_{1}^{s})$,
$\,i=1,\ldots,n$, for some $y_{1}^{s}\in G^{s},$
$\,i_1,\ldots,i_s\in\{1,\ldots,n\}$, then $\;h_1\sqsubset
g_1\opp{i_{1}}{i_{s}}y_{1}^{s}\,\vee\, h_2\sqsubset
g_1\opp{i_{1}}{i_{s}}y_{1}^{s}\;$ and $\;h_1\sqsubset
g_2\opp{i_{1}}{i_{s}}y_{1}^{s}\,\vee\,h_2\sqsubset
g_2\opp{i_{1}}{i_{s}}y_{1}^{s}$, whence, by $\,h_1\top\,h_2$ and
(\ref{15}), we obtain
$\,g_1\opp{i_{1}}{i_{s}}y_{1}^{s}\,\top\,g_2\opp{i_{1}}{i_{s}}y_{1}^{s}$.
This gives $g_1\top\,g_2$ because $\gamma$ is $l$-cancellative. In
this way we have proved that in any case $\gamma_P\subset\gamma$.

Conversely, let $(g_1,g_2)\in\gamma$. Since $\chi$ is reflexive,
$g_1\sqsubset g_1$ and $g_2\sqsubset g_2$, whence $g_1\sqsubset
g_1\,\vee\, g_2\sqsubset g_1\,$ and $\,g_1\sqsubset
g_2\,\vee\,g_2\sqsubset g_2$. Consequently, $e_1^n\in\pr
P_{(g_{1},\,g_{2})}(g_1)$ and $e_1^n\in\pr
P_{(g_{1},\,g_{2})}(g_2)$. Thus
$(g_1,g_2)\in\gamma_{(g_1,g_2)}\subset\gamma_P$, i.e.
$\gamma\subset\gamma_P$. So, $\gamma=\gamma_P$.

This completes the proof of the theorem. \hfill$\Box$

\bigskip

\noindent{\bf Problem 1.} {\it Find the necessary and sufficient
conditions under which the triplet $(\chi,\gamma,\pi)$ of binary
relations will be faithful projection representable for a
representable Menger $(2,n)$-semigroup.}

\medskip

Deleting from Theorem~\ref{T1} the equality
$\pi=\chi\cap\chi^{-1}$ we obtain the necessary and sufficient
conditions under which the pair $(\chi,\gamma)$ of binary
relations is projection representable for a representable Menger
$(2,n)$-semigroup. Furthermore, all parts of the proof of this
theorem connected with these two relations are valid. So, we have
the following

\begin{Theorem}\label{T1a}
A pair $(\chi,\pi)$ of binary relations on $G$ is projection
representable for a representable Menger $(2,n)$-semigroup
$\mathcal{G}$ if and only if $\chi$ is an $l$-regular and
$v$-negative quasi-order, $\gamma$ is an $l$-cancellative
$0$-quasi-equivalence and the implication $(\ref{15})$ is
satisfied.
\end{Theorem}

\bigskip

\noindent{\bf Problem 2.} {\it Find the necessary and sufficient
conditions under which the pair $(\chi,\gamma)$ of binary
relations will be faithful projection representable for a
representable Menger $(2,n)$-semigroup.}

\bigskip

Let $\mathcal{G}=(G;o,\op{1\;},\ldots,\op{n})$ be a representable
Menger $(2,n)$-semigroup. Let us consider on $G$ the set $T_n(G)$
of mappings $\;t:x\mapsto t(x)$ defined as follows:
\begin{itemize}
  \item[(a)] $x\in T_n(G)$, \ i.e. $T_n(G)$ contains the identity
transformation of $G$,
  \item[(b)] if $i\in\{1,\ldots,n\}$, $\,a,b_1,\ldots,b_{i-1},b_{i+1},\ldots,b_n\in
  G$ and $t(x)\in T_n$, then $a[b_{1}^{i-1}t(x)\,b_{i+1}^{n}]\in
  T_n$,
\item[(c)] $T_n$ contains those and only those mappings which are defined by (a) and
  (b).
\end{itemize}
Let us consider on $G$ two binary relations $\delta_1$ and
$\delta_2$ defined in the following way:

\begin{enumerate}
  \item $(g_1,g_2)\in\delta_1\,\longleftrightarrow\,g_1=t(g_2)$ for some $t\in T_n$,
  \item $(g_1,g_2)\in\delta_2\,\longleftrightarrow
  \left\{\begin{array}{l} g_1=(x\opp{i_1}{i_s}y_1^s)[\bar{z}] \ \ {\rm and}
  \ \ g_2=\mu_i(\opp{i_1}{i_s}y_1^s)[\bar{z}] \ \ {\rm for\;some }\\
  x\in G, \ y_1^s\in G^s, \ \bar{z}\in G^n, \
  i,i_1,\ldots,i_s\in\{1,\ldots,n\},\\[4pt]
  {\rm where\; the\; symbol } \ [\bar{z}] \ {\rm can\; be\; empty}.
  \end{array}\right.$
\end{enumerate}
It is not difficult to see that $\delta_1$ and $\delta_2$ are
$l$-regular relations, additionally $\delta_1$ is a quasi-order.
Moreover, {\it a binary relation $\rho\subset G\times G$ is
$v$-negative if and only if it contains $\delta_1$ and
$\delta_2$}.

Let $\pi$ be an $l$-regular equivalence on a representable Menger
$(2,n)$-semigroup $\mathcal{G}$. Denote by $\chi(\pi)$ the binary
relation $f_t(f_R(\delta_2)\circ\delta_1\circ\pi)$, where $f_R$
and $f_t$ are respectively reflexive and transitive closure
operations (cf. \cite{Rig}), and $\circ$ is a composition of
relations, \footnote{ Remind that
$\sigma\circ\rho=\{(a,c)\,|\,(\exists
b)(a,b)\in\rho\,\wedge\,(b,c)\in\sigma\}$,
$\,f_R(\rho)=\rho\cup\bigtriangleup_A$,
$\,f_t(\rho)=\bigcup\limits_{n=1}^{\infty}\rho^n$, where
$\rho^n=\underbrace{\rho\circ\rho\circ\ldots\circ\rho}_n$,
$\;\rho,\sigma$~--- binary relations on $A$ and
$\triangle_A=\{(a,a)\,|\,a\in A\}$.} i.e.
\begin{equation}\label{32}
\chi(\pi)=f_t(f_R(\delta_2)\circ\delta_1\circ\pi)=
\bigcup\limits_{n=1}^{\infty}\left((\delta_2\cup\triangle_G)\circ
\delta_1\circ\pi\right)^{n}.
\end{equation}
Since $\pi$, $\delta_1$ and $f_R(\delta_2)$ are reflexive
$l$-regular relations, $\chi(\pi)$ is an $l$-regular quasi-order
containing $\pi$, $\delta_1$ and $\delta_2$. So, $\chi(\pi)$ is a
$v$-negative quasi-order.

\begin{Proposition}\label{P4}
$\chi(\pi)$ is the least $l$-regular and $v$-negative quasi-order
containing $\pi$.
\end{Proposition}

\begin{proof} Let $\chi$ be an arbitrary
$l$-regular and $v$-negative quasi-order containing $\pi$. Then
$\delta_1\subset\chi$ and $\delta_2\subset\chi$, because $\chi$ is
$v$-negative. Thus, $\pi\subset\chi$, $\delta_1\subset\chi$ and
$f_R(\delta_2)\subset\chi$, whence
$f_R(\delta_2)\circ\delta_1\circ\pi\subset\chi^3\subset\chi$. From
this, applying the transitivity of $\chi$, we obtain
$(f_R(\delta_2)\circ\delta_1\circ\pi)^n\subset\chi^n\subset\chi$
for every natural $n$. Therefore
$\bigcup\limits_{n=1}^{\infty}\left((\delta_2\cup\triangle_G)\circ\delta_1\circ\pi\right)^n
\subset\chi$, i.e. $\chi(\pi)\subset\chi$.
 \end{proof}

\begin{Theorem}\label{T2}
A pair $(\gamma,\pi)$ of binary relations on a representable
Menger $(2,n)$-semigroup $\mathcal{G}$ is projection representable
if and only if

$(a)$ \ $\gamma$ is an $l$-cancellative $0$-quasi-equivalence,

$(b)$ \ $\pi$ is an $l$-regular equivalence such that
$\;\chi(\pi)\cap(\chi(\pi))^{-1}\subset\pi$,

$(c)$ \ the following condition
\begin{equation}\label{33}
  h_1\top\,h_2\,\wedge\, h_1\sqsubset_{\,\pi}g_1\,\wedge\,
  h_2\sqsubset_{\,\pi}g_2\longrightarrow g_1\top\,g_2,
\end{equation}
where $h\sqsubset_{\,\pi}g$ means $(h,g)\in\chi(\pi)$, is
satisfied for all $g_1,g_2,h_1,h_2\in G$.
\end{Theorem}

\begin{proof} Let $P$ be such representation on a representable
Menger $(2,n)$-semigroup $\mathcal{G}$ for which $\gamma=\gamma_P$
and $\pi=\pi_P$. Then, by Proposition~\ref{P3}, we have
$\chi(\pi)\subset\chi_P$, whence
$\chi(\pi)\cap(\chi(\pi))^{-1}\subset\chi_P\cap\chi_{P}^{-1}=\pi_P=\pi$.

Assume now that the premise of (\ref{33}) is satisfied. Then
$(h_1,h_2)\in\gamma$, $\,(h_1,g_1)\in\chi(\pi)$ and
$(h_2,g_2)\in\chi(\pi)$. Consequently, $(h_1,h_2)\in\gamma_P$,
$\,(h_1,g_1)\in\chi_P$ and $(h_2,g_2)\in\chi_P$, \ i.e. $\pr
P(h_1)\cap\pr P(h_2)\neq\varnothing$, $\,\pr P(h_1)\subset\pr
P(g_1)$ and $\pr P(h_2)\subset\pr P(g_2)$, whence $\pr
P(g_1)\cap\pr P(g_2)\neq\varnothing$. So,
$(g_1,g_2)\in\gamma_P=\gamma$, which means that the condition
(\ref{33}) is valid. The necessity is proved.

To prove the sufficiency, assume that the pair $(\gamma,\pi)$ of
binary relations satisfies all the conditions of the theorem and
consider the triplet $(\chi(\pi),\gamma,\pi)$. Then
$\pi=\pi^{-1}\subset(\chi(\pi))^{-1}$, because
$\pi\subset\chi(\pi)$. Therefore
$\pi\subset\chi(\pi)\cap(\chi(\pi))^{-1}$, which, together with
the condition $(b)$, gives $\pi=\chi(\pi)\cap(\chi(\pi))^{-1}$.
This means that the triplet $(\chi(\pi),\gamma,\pi)$ satisfies all
the conditions of Theorem~\ref{T1}. So, $(\chi(\pi),\gamma,\pi)$,
and in the consequence, $(\gamma,\pi)$ is projection
representable. The sufficiency is proved.
\end{proof}

\noindent{\bf Problem 3.} {\it Find the necessary and sufficient
conditions under which the pair $(\gamma,\pi)$ of binary relations
will be faithful projection representable.}

\bigskip

Applying the method of mathematical induction to (\ref{32}) we can
prove the following proposition.

\begin{Proposition}\label{P5}
The condition $(g_1,g_2)\in\chi(\pi)$, where $g_1,g_2\in G$, means
 that the system of conditions
\begin{equation}\label{34}
  \left.
\begin{array}{c}
  g_1=x_0\,\wedge\, g_2=x_n,\\[4pt]
  \bigwedge\limits_{i=0}^{n-1}\left(\left(
\begin{array}{l}
  x_i\equiv t_i((y_i\opp{k_{1_i}}{k_{s_i}}z_{1_i}^{s_i})[\bar{w}_i]), \\[4pt]
  x_{i+1}=\mu_{k_i}(\opp{k_{1_i}}{k_{s_i}}z_{1_i}^{s_i})[\bar{w}_i]
\end{array}
  \right)\,\vee\, x_i\equiv t_i(x_{i+1})\right)
\end{array}
  \right\}
\end{equation}
is valid for some $\,n\in\mathbb{N}$, $\,x_i,y_i,z_i\in G$,
$\,\bar{w}_i\in G^n$, $\,t_i\in T_n$, $\,k_i\in\{1,\ldots,n\}$.
\qed
\end{Proposition}

\medskip

In the sequel the formula
$$
\bigwedge\limits_{i=m}^{n}\left(\left(
\begin{array}{l}
  x_i\equiv t_i((y_i\opp{k_{1_i}}{k_{s_i}}z_{1_i}^{s_i})[\bar{w}_i]), \\[4pt]
  x_{i+1}=\mu_{k_i}(\opp{k_{1_i}}{k_{s_i}}z_{1_i}^{s_i})[\bar{w}_i]
\end{array}
  \right)\vee x_i\equiv t_i(x_{i+1})\right)
$$
will be denoted by $\frak{M}(m,n)$.

The inclusion $\,\chi(\pi)\cap(\chi(\pi))^{-1}\subset\pi\,$ means
that for all $g_1,g_2\in G$ we have
\[
(g_1,g_2)\in\chi(\pi)\,\wedge\,(g_2,g_1)\in\chi(\pi)\longrightarrow
g_1\equiv g_2,
\]
which, according to Proposition~\ref{P5}, can be written as the
system of conditions $(A_{n,m})_{n,m\in\mathbb{N}}$, where
\[
A_{n,m}:\ \frak{M}(0,n-1)\,\wedge\,\frak{M}(n+1,n+m)\,\wedge\,
x_0=x_{n+m}\longrightarrow x_0\equiv x_n.
\]
The system $(A_{n,m})_{n,m\in\mathbb{N}}$ is equivalent to the
system $(A_{n})_{n\in\mathbb{N}}$, where
\[
A_n:\ \frak{M}(0,n-1)\,\wedge\, x_0=x_n\longrightarrow x_0\equiv
x_1.\]

Consider now the implication (\ref{33}). According to (\ref{34})
the condition $(h_1,g_1)\in\chi(\pi)$ means that
\begin{equation}\label{35}
  h_1=x_0\,\wedge\,\frak{M}(0,n-1)\,\wedge\, x_n=g_1
\end{equation}
for some $x_i,y_i,z_{k_i},t_i,k_i,\bar{w}_i$. Similarly, the
condition $(h_2,g_2)\in\chi(\pi)$ means that
\begin{equation}\label{36}
  h_2=x_{n+1}\,\wedge\,\frak{M}(n+1,n+m)\,\wedge\, x_{n+m+1}=g_2
\end{equation}
for some $x_i,y_i,z_{k_i},t_i,k_i,\bar{w}_i$. So, (\ref{33}) can
be written as the system $(B_{n,m})_{n,m\in\mathbb{N}}$ of
conditions
\[
B_{n,m}:\
x_0\top\,x_{n+1}\,\wedge\,\frak{M}(0,n-1)\,\wedge\,\frak{M}(n+1,n+m)
\longrightarrow x_n\top\,x_{n+m+1}.
\]

In this way we have proved

\begin{Theorem}\label{T3}
A pair $(\gamma,\pi)$ of binary relations on a representable
Menger $(2,n)$-semigroup $\mathcal{G}$ is projection representable
if and only if

$(a)$ \ $\gamma$ is an $l$-cancellative $0$-quasi-equivalence,

$(b)$ \ $\pi$ is an $l$-regular equivalence,

$(c)$ \ the systems of conditions $(A_n)_{n\in\mathbb{N}}$ and
$(B_{n,m})_{n,m\in\mathbb{N}}$ are satisfied. \qed
\end{Theorem}

\begin{Theorem}\label{T4}
A pair $(\chi,\pi)$ of binary relations is $($faithful$)$
projection represen\-table for a representable Menger
$(2,n)$-semigroup $\mathcal{G}$ if and only if $\chi$ is an
$l$-regular and $v$-negative quasi-order such that
$\pi=\chi\cap\chi^{-1}$.
\end{Theorem}

\begin{proof}
The necessity of these conditions follows from the proof of
Theorem~\ref{T1}. To prove their sufficiency, for every element
$g\in G$ we define an $n$-place function $P_a(g):
\frak{A}^{*}\rightarrow G$, where $a\in G$, putting
\begin{equation}\label{37}
P_a(g)(x_1^n)=\left\{
\begin{array}{ll}
 g[x_1^n] &\mbox{if } a\sqsubset g[x_1^n] \ \mbox{ and } x_1^n\in G^n, \\[4pt]
 g & \mbox{if } a\sqsubset g \ \mbox{ and } \ x_1^n=e_1^n, \\[4pt]
g\opp{i_1}{i_s}y_1^s &\mbox{if } a\sqsubset g\opp{i_1}{i_s}y_1^s \
\mbox{ and } \ x_i=\mu_i^*(\opp{i_1}{i_s}y_1^s), \\
 & i=1,\ldots,n,\mbox{ for\; some }y_1^s\in G^s, \\[3pt]
 & i_1,\ldots,i_s\in\{1,\ldots,n\}.
\end{array}
\right.
\end{equation}
Since, for $h_1=h_2=a\in G,$ the function $P_{(h_1,h_2)}(g)$
defined by (\ref{16}) coincides with the function $P_a(g)$, from
Propositions~\ref{P1}~--~\ref{P3} it follows that the mapping
$P_a:g\mapsto P_a(g)$ is a representation of $\mathcal{G}$ by
$n$-place functions. Further, analogously as in the proof of
Theorem~\ref{T1}, we can prove that $P_0=\sum\limits_{a\in G}P_a$
is a representation of $\mathcal{G}$ for which $\chi=\chi_{P_0}$
and $\pi=\pi_{P_0}$. So, the pair $(\chi,\pi)$ is projection
representable for $\mathcal{G}$.

Let us show that $(\chi,\pi)$ is faithful projection
representable. In \cite{Sok} it is proved that each representable
Menger $(2,n)$-semigroup has a faithful representation by
$n$-place functions. Let $\Lambda$ be such representation. Then
obviously $\chi_{\Lambda}=G\times G$ and $\pi_{\Lambda}=G\times
G$. Consider the representation $P=\Lambda+P_0$. Since $\Lambda$
is a faithful representation, $P$ is also faithful. Moreover
$\chi_P=\chi_{\Lambda}\cap\chi_{P_0}=G\times G\cap\chi=\chi$ and
$\pi_P=\pi_{\Lambda}\cap\pi_{P_0}=G\times G\cap\pi=\pi$. So,
$(\chi,\pi)$ is faithful projection representable for
$\mathcal{G}$.
 \end{proof}

In the same manner, using the construction (\ref{37}), we can
prove the following theorem.
\begin{Theorem}\label{T5}
A binary relation $\chi$ is $($faithful$)$ projection
representable for a representable Menger $(2,n)$-semigroup if and
only if it is an $l$-regular, $v$-negative quasi-order. \qed
\end{Theorem}
\begin{Theorem}\label{T6}
A binary relation $\pi$ is $($faithful$)$ projection representable
for a representable Menger $(2,n)$-semigroup if and only if it is
an $l$-regular equivalence such that
$\;\chi(\pi)\cap(\chi(\pi))^{-1}\subset\pi$.
\end{Theorem}

\begin{proof} Consider the pair $(\chi(\pi),\pi)$ of binary relations, where
$\chi(\pi)$ is defined by (\ref{32}). In a similar way, as in the
proof of Theorem~\ref{T2}, we can prove that this pair satisfies
all the conditions of Theorem~\ref{T4}, whence we conclude the
validity of Theorem~\ref{T6}.
 \end{proof}

Since, as it was showed above, the inclusion
$\,\chi(\pi)\cap(\chi(\pi))^{-1}\subset\pi$ is equivalent to the
system of conditions $(A_n)_{n\in\mathbb{N}}$, the last theorem
can be rewritten in the form:
\begin{Theorem}\label{T7}
A binary relation $\pi$ is $($faithful$)$ projection representable
for a representable Menger $(2,n)$-semigroup if and only if it is
an $l$-regular equivalence and the system of conditions
$(A_n)_{n\in\mathbb{N}}$ is satisfied. \qed
\end{Theorem}

Consider on a Menger $(2,n)$-semigroup $\mathcal{G}$ the binary
relation $\chi_0$ defined in the following way:
\begin{equation}\label{38}
\chi_0=f_t(f_R(\delta_2)\circ\delta_1)=
\bigcup\limits_{n=1}^{\infty}
\left((\delta_2\cup\bigtriangleup_G)\circ\delta_1\right)^{n},
\end{equation}
where $f_t$ and $f_R$ are reflexive and transitive closure
operations.
\begin{Proposition}\label{P6}
$\chi_0$ is the least $l$-regular and $v$-negative quasi-order on
$\mathcal{G}$.
\end{Proposition}

The proof of this proposition is analogous to the proof of
Proposition~\ref{P3}.

\begin{Theorem}\label{T8}
A binary relation $\gamma$ is projection representable for a
representable Menger $(2,n)$-semigroup if and only if it is an
$l$-cancellative $0$-quasi-equivalence and the following
implication
\begin{equation}\label{39}
h_1\top\,h_2\,\wedge\, h_1\sqsubset_{\,0} g_1\,\wedge\,
h_2\sqsubset_{\,0} g_2\;\longrightarrow g_1\top\,g_2
\end{equation}
is satisfied for all $\,h_1,h_2,g_1,g_2\in G$, where
$h\sqsubset_{\,0} g$ means $(h,g)\in\chi_0$.
\end{Theorem}

\begin{proof} The necessity of (\ref{39}) can be proved
analogous as the necessity of (\ref{33}) in the proof of
Theorem~\ref{T2}. To prove the sufficiency we consider the pair
$(\chi_0,\gamma)$. By Proposition~\ref{P6}, this pair satisfies
all demands of Theorem~\ref{T1a}, whence we conclude the validity
of Theorem~\ref{T8}.
 \end{proof}

\noindent{\bf Problem 4.} {\it Find the necessary and sufficient
conditions under which $\gamma$ will be faithful projection
representable.}

\bigskip

Basing on the formula (\ref{38}) we can prove the following
proposition:
\begin{Proposition}\label{P7}
From $(g_1,g_2)\in \chi_0$, where $g_1,g_2\in G$,  it follows that
the system of conditions
\begin{equation*}\label{40}
  g_1=x_0\,\wedge\, g_2=x_n\,\wedge\,
  \bigwedge\limits_{i=0}^{n-1}\left(\left(
\begin{array}{l}
  x_i=t_i((y_i\opp{k_{1_i}}{k_{s_i}}z_{1_i}^{s_i})[\bar{w}_i]), \\[4pt]
  x_{i+1}=\mu_{k_i}(\opp{k_{1_i}}{k_{s_i}}z_{1_i}^{s_i})[\bar{w}_i]
\end{array}
  \right)\vee\, x_i=t_i(x_{i+1})\right)
\end{equation*}
is valid for $n\in\mathbb{N}$, $\,x_i,y_i,z_i\in G,$
$\,\bar{w}_i\in G^n,$ $\,t_i\in T_n$, $\,k_i\in\{1,\ldots,n\}$.
\qed
\end{Proposition}

Denoting by $\frak{N}(m,n)$ the formula
\[
\bigwedge\limits_{i=m}^{n}\left(\left(
\begin{array}{l}
  x_i= t_i((y_i\opp{k_{1_i}}{k_{s_i}}z_{1_i}^{s_i})[\bar{w}_i]), \\[4pt]
  x_{i+1}=\mu_{k_i}(\opp{k_{1_i}}{k_{s_i}}z_{1_i}^{s_i})[\bar{w}_i]
\end{array}
  \right)\vee\, x_i=t_i(x_{i+1})\right),
\]
and using the same argumentation as in the proof of
Theorem~\ref{T3}, we can prove that the implication (\ref{39}) is
equivalent to the system of conditions
$(C_{n,m})_{n,m\in\mathbb{N}}$, where
\[
C_{n,m}:\
x_0\top\,x_{n+1}\,\wedge\,\frak{N}(0,n-1)\,\wedge\,\frak{N}(n+1,n+m)\,
\longrightarrow x_0\top\,x_{n+m+1}.
 \]

So, the following theorem is true:
\begin{Theorem}\label{T9}
A binary relation $\gamma$ is projection representable for a
representable Menger $(2,n)$-semigroup if and only if it is an
$l$-cancellative $0$-quasi-equivalence and the system of
conditions $(C_{n,m})_{n,m\in\mathbb{N}}$ is satisfied. \qed
\end{Theorem}

\section{Projection representable relations on
$(2,n)$-semi\-groups}

Let $\chi$, $\gamma$ and $\pi$ be three binary relations on a
$(2,n)$-semigroup $(G;\op{1\,},\ldots,\op{n})$. Similarly as in
the case of Menger $(2,n)$-semigroups we say that the triplet
$(\chi,\gamma,\pi)$ is ({\it faithful}) \textit{projection
representable for a $(2,n)$-semigroup
$(G;\op{1\,},\ldots,\op{n})$}, if there exists such (faithful)
representation $P$ of $(G;\op{1\,},\ldots,\op{n})$ by $n$-place
functions for which $\chi=\chi_{P}$, $\gamma=\gamma_{P}$ and
$\pi=\pi_{P}$. Analogously we define the projection representable
pairs and separate relations.

\medskip

It is not difficult to verify that our Theorem~\ref{T1} formulated
for representable Menger $(2,n)$-semigroup is also valid for
representable $(2,n)$-semigroups. The proof of this version of
Theorem~\ref{T1} is analogous to the proof of the previous
version, but in the proof of the sufficiency instead the
representation $P$ we must consider the representation
$P^{\,\bullet}$, which is the sum of the family of representations
$(P^{\,\bullet}_{(h_1,\,h_2)})_{(h_1,\,h_2)\in\gamma}$, where for
every $g\in G$ \ $P^{\,\bullet}_{(h_1,\,h_2)}(g):\
\frak{A}^*_0\rightarrow G$
($\frak{A}^*_0=\frak{A}_0\cup\{(e_1,\ldots,e_n)\}$ see
page~\pageref{Tr**}) is a partial $n$-place function such that
\[
x_{1}^{n}\in\pr
P^{\,\bullet}_{(h_{1},\,h_{2})}(g)\,\longleftrightarrow
 \left\{
\begin{array}{ll}
h_{1}\sqsubset g\,\vee\,h_{2}\sqsubset g & \mbox{if } x_{1}^{n}=e_{1}^{n},\\[4pt]
h_{1}\sqsubset g\opp{i_{1}}{i_{s}}y_{1}^{s}\,\vee\, h_{2}\sqsubset
g\opp{i_{1}}{i_{s}}y_{1}^{s} & \mbox{if }\
  x_{i}=\mu_{i}^{*}(\opp{i_{1}}{i_{s}}y_{1}^{s}),\\
  & i=1,\ldots,n,\mbox{ for}  \\
  & \mbox{some } y_{1}^{s}\in G^{s} \mbox{ and } \\
  & i_{1}\ldots,i_{s}\in\{1,\ldots,n\}
\end{array}
  \right.
\]
and
 \[
P^{\,\bullet}_{(h_{1},\,h_{2})}(g)(x_{1}^{n})= \left\{
\begin{array}{ll}
    g & \mbox{if } x_{1}^{n}=e_{1}^{n}, \\[4pt]
  g\opp{i_{1}}{i_{s}}y_{1}^{s} \ \ & \mbox{if }\
  x_{i}=\mu_{i}^{*}(\opp{i_{1}}{i_{s}}y_{1}^{s}),\\[1pt]
  & i=1,\ldots,n,\mbox{ for\;some }y_{1}^{s}\in G^{s}\\ & \mbox{and
  }\ i_{1}\ldots,i_{s}\in\{1,\ldots,n\}.
\end{array}
  \right.
\]

\medskip

Also Theorem~\ref{T1a} is valid for $(2,n)$-semigroups. Moreover,
problems analogous to Problem 1 and Problem 2 can be posed for
$(2,n)$-semigroups, too.

\medskip

Theorem~\ref{T2} will be valid for $(2,n)$-semigroups if we
replace the relation $\chi(\pi)$ by the relation
\begin{equation}\label{62}
\chi^{\,\bullet}(\pi)=f_t(f_R(\delta_2)\circ\pi)=
\bigcup\limits_{n=1}^{\infty}\left((\delta_2\cup\triangle_G)\circ\pi\right)^{n},
\end{equation}
i.e. if we delete $\delta_1$ from the formula (\ref{32}).

\medskip

Proposition~\ref{P5} for $(2,n)$-semigroups has the following
form:
\begin{Proposition}\label{P8}
The condition $(g_1,g_2)\in\chi^{\,\bullet}(\pi)$, where
$g_1,g_2\in G,$ means that the system of conditions
\begin{equation*}
  g_1=x_0\,\wedge\, g_2=x_n\,\wedge\,
  \bigwedge\limits_{i=0}^{n-1}\left(\left(
\begin{array}{c}
  x_i\equiv y_i\opp{k_{1_i}}{k_{s_i}}z_{1_i}^{s_i}, \\[4pt]
  x_{i+1}=\mu_{k_i}(\opp{k_{1_i}}{k_{s_i}}z_{1_i}^{s_i})
\end{array}
  \right)\vee\, x_i\equiv x_{i+1}\right)
\end{equation*}
is valid for some $n\in\mathbb{N}$, $\,x_i,y_i,z_i\in G,$
$\,k_i\in\{1,\ldots,n\}$. \qed
\end{Proposition}

Denoting by $\frak{X}(m,n)$ the formula
 \[
\bigwedge\limits_{i=m}^{n}\left(\left(
\begin{array}{l}
  x_i\equiv y_i\opp{k_{1_i}}{k_{s_i}}z_{1_i}^{s_i}, \\[4pt]
  x_{i+1}=\mu_{k_i}(\opp{k_{1_i}}{k_{s_i}}z_{1_i}^{s_i})
\end{array}
  \right)\vee\, x_i\equiv x_{i+1}\right)
\]
and using the same argumentation as in the proof of
Theorem~\ref{T3}, we can prove
\begin{Theorem}\label{T11}
A pair $(\gamma,\pi)$ of binary relations on a representable
$(2,n)$-semi\-group is projection representable if and only if
$\gamma$ is an $l$-cancellative $0$-quasi-equivalence, $\pi$ is an
$l$-regular equivalence, and the systems of conditions
$A^{\,\bullet}_n$ and $B^{\,\bullet}_{n,m}$, where
\begin{eqnarray*}
&&A^{\,\bullet}_n:\  \frak{X}(0,n-1)\,\wedge\,
x_0=x_n\longrightarrow
x_0\equiv x_1,\\[4pt]
&&B^{\,\bullet}_{n,m}:\
x_0\top\,x_{n+1}\,\wedge\,\frak{X}(0,n-1)\,\wedge\,\frak{X}(n+1,n+m)\longrightarrow
x_n\top\,x_{n+m+1}
\end{eqnarray*}
are satisfied. \qed
\end{Theorem}

Theorem~\ref{T4} is valid for $(2,n)$-semigroups too, but in the
proof, the representation $P_a$ defined by (\ref{37}), must be
replaced by the representation $P^{\,\bullet}_a$, where
\begin{equation*}
P^{\,\bullet}_a(g)(x_1^n)=\left\{
\begin{array}{ll}
 g & \mbox{if } a\sqsubset g\,\mbox{ and }\,x_1=e_1^n, \\[4pt]
 g\opp{i_1}{i_s}y_1^s & \mbox{if } a\sqsubset g\opp{i_1}{i_s}y_1^s\,\mbox{ and }\,x_i=
 \mu_i^*(\opp{i_1}{i_s}y_1^s), \\
 & i=1,\ldots,n,\mbox{ for some }\,y_1^s\in G^s, \\[2pt]
    & \mbox{and }\,i_1,\ldots,i_s\in\{1,\ldots,n\}.
\end{array}
\right.
\end{equation*}

For $(2,n)$-semigroups Theorem~\ref{T5} has the same form as for
Menger $(2,n)$-semigroup, in Theorem~\ref{T6} the relation
$\chi(\pi)$ must be replaced by $\chi^{\,\bullet}(\pi)$, and in
Theorem~\ref{T7} instead of $A_n$ we must use $A^{\,\bullet}_n$.

Further, using the same argumentation as in the proof of
Proposition~\ref{P4} we can prove that the relation
\begin{equation*}
\chi^{\,\bullet}_0=f_t(f_R(\delta_2))=
\bigcup\limits_{n=1}^{\infty}\left(\delta_2\cup\bigtriangleup_G\right)^{n},
\end{equation*}
where $f_t$ and $f_R$ are reflexive and transitive closure
operations, is the least $l$-regular and $v$-negative quasi-order
on a given $(2,n)$-semigroup. Using this relation, we can prove
the analog of Theorem~\ref{T9} for $(2,n)$-semigroups. The analog
of Problem 4 can be posed too.

Proposition~\ref{P7} for $(2,n)$-semigroups has the following
form:
\begin{Proposition}\label{P9} The condition
$(g_1,g_2)\in\chi^{\,\bullet}_0$, where $g_1,g_2\in G$, means that
the system of conditions
\begin{equation*}
  g_1=x_0\,\wedge\, g_2=x_n\,\wedge\,
  \bigwedge\limits_{i=0}^{n-1}\left(\left(
\begin{array}{c}
  x_i=y_i\opp{k_{1_i}}{k_{s_i}}z_{1_i}^{s_i}, \\[4pt]
  x_{i+1}=\mu_{k_i}(\opp{k_{1_i}}{k_{s_i}}z_{1_i}^{s_i})
\end{array}
  \right)\vee x_i=x_{i+1}\right)
\end{equation*}
is valid for $n\in\mathbb{N}$, $\,x_i,y_i,z_i\in G$. \qed
\end{Proposition}

Further, denoting by $\frak{B}(m,n)$ the formula
$$
\bigwedge\limits_{i=m}^{n}\left(\left(
\begin{array}{l}
  x_i=y_i\opp{k_{1_i}}{k_{s_i}}z_{1_i}^{s_i}, \\[4pt]
  x_{i+1}=\mu_{k_i}(\opp{k_{1_i}}{k_{s_i}}z_{1_i}^{s_i})
\end{array}
  \right)\vee x_i=x_{i+1}\right)
$$
and using the same argumentation as in the proof of
Theorem~\ref{T9}, we can prove
\begin{Theorem}\label{T12} A binary
relation $\gamma$ is projection representable for a representable
$(2,n)$-semigroup if and only if it is an $l$-cancellative
$0$-quasi-equivalence and the system of conditions
$(C^{\,\bullet}_{n,m})_{n,m\in\mathbb{N}}$, where
\[
C^{\,\bullet}_{n,m}:\
x_0\top\,x_{n+1}\,\wedge\,\frak{B}(0,n-1)\,\wedge\,\frak{B}(n+1,n+m)\,
\longrightarrow x_0\top\,x_{n+m+1}
\]
is satisfied. \qed
\end{Theorem}

\medskip\noindent
{\sc Wies{\l}aw A. Dudek}\\
Institute of Mathematics, Technical University, 50-370
Wroc{\l}aw, Poland\\
E-mail: dudek@im.pwr.wroc.pl\\[10pt]
{\sc Valentin S. Trokhimenko}\\
Department of Mathematics, Pedagogical University,
21100 Vinnitsa, Ukraine\\
E-mail: vtrokhim@sovamua.com

\end{document}